# Safe control of thermostatically controlled loads with installed timers for demand side management

Nishant Mehta, Nikolai. A. Sinitsyn, Scott Backhaus, Bernard. C. Lesieutre

*Abstract*–We develop safe protocols for thermostatically controlled loads (TCLs) to provide power pulses to the grid without a subsequent oscillatory response. Such pulses can alleviate power fluctuations by intermittent resources and maintain balance between generation and demand. Building on prior work, we introduce timers to endpoint TCL control enabling better shaping of power pulses.

*Index Terms*-Ancillary services, Demand response, Demand side management, Generation-load balance, Hysteresis-based control, Load control, Load modeling, Power demand, Renewable energy, Thermostatically controlled loads.

## I. INTRODUCTION

When generators are constrained by their ramp rates, it is difficult to adjust the generation to meet sudden peaks in demand without a significant delay [1]. Adding more generation capacity may not be an optimal solution to this problem because of additional costs and emissions. It has been proposed previously that thermostatically controlled loads (TCLs), such as air-conditioners, refrigerators and water heaters, can provide ancillary services to balance generation and demand [2–16].

Many TCLs turn between ON and OFF states depending on their ambient conditions [1], [8], [14], [15], [17–34]. TCLs are considered as appropriate candidates to provide ancillary services for two main reasons – (i) around 50% of electricity in the United States is consumed by TCLs, (ii) TCLs have inherent ability to store thermal energy in the space they are managing [25]. According to US Energy Information Administration [35], American households that have air-conditioners have increased their number from 68% in 1993 to 87% in 2009 with 65.1% of all occupied homes in US having central air-conditioning units. In China alone, over 50 million of air conditioners are sold annually [36]. As TCLs are able to respond to control signals faster than spinning reserves [9], [16], [37–39], they can help to stabilize the power balance in the grid by offsetting sharp power fluctuations at time scales of minutes.

The potential of TCLs to assist in controlling the power grid by varying their temperature set points has been discussed in a number of publications [17], [22], [24], [25], [28], [30]. A major challenge to achieving control over a large ensemble of TCLs is the wide distribution of parameters such as thermal resistance, thermal capacitance, set-point, ambient temperature, duty cycle and the power consumed when an air-conditioner is in the ON state. It was discovered in numerical studies that, because of the lack of complete information, it is hard to avoid unwanted power oscillations lasting for several hours [25], [28], due to synchronization of TCLs by the external control signals.

The aim of this paper is to explore approaches to constructing safe protocols such as adding timers at the endpoint TCL temperature controllers in order to generate on-demand power outputs by TCL ensembles while having minimal information about the aggregate TCL ensemble parameters. Adding memory and instructions to the endpoint TCL temperature controllers would enable them to autonomously switch between ON/OFF states without waiting for a specific external signal. Specifically, we will (i) propose a safe protocol not discussed in the previous papers [1], [20], that generates a power pulse using a timer, (ii) we perform simulations to demonstrate that timer-based safe protocols can be used to successfully eliminate strong and short power fluctuations in the power grid and (iii) we explore the application of a safe protocol beyond the simple framework of a steady state equilibrated ensemble of TCLs.

The rest of the paper is organized as follows. Section II describes the model of a TCL that controls the conditioned space temperature. Section III demonstrates the *'Unsafe'* set point change method to control the aggregated power. Section IV describes *'Safe Protocols'* to generate power pulses. Section V demonstrates the versatility of safe protocols in offsetting fast time scale fluctuations. Section VI explains the role of heterogeneity where ambient temperature is varying throughout the day, and Section VII concludes the paper.

## II. MODEL OF TCL DYNAMICS

Let the TCL be controlling the room temperature $\theta(t)$. The differential equation that describes the evolution of $\theta(t)$ within the dead band limits $[\theta_-, \theta_+]$, where $\theta_-$ and $\theta_+$ being the lower and upper dead band limits respectively, are as follows:

Nishant Mehta and Bernard Lesieutre are with the Department of Electrical and Computer Engineering at the University of Wisconsin-Madison, Madison, WI-53706 USA (e-mail: nishantmehta1289@gmail.com, lesieutre@engr.wisc.edu).

Nikolai Sinsitsyn and Scott Backhaus are with the Los Alamos National Laboratory (LANL), Los Alamos, NM-87545 USA (e-mail: nsinitsyn@lanl.gov, backhaus@lanl.gov).

$$\dot{\theta} = \begin{cases} \dfrac{-1}{CR}\left(\theta - \theta_{amb} + PR + \xi(t)\right) & \text{ON state} \\[2mm] \dfrac{-1}{CR}\left(\theta - \theta_{amb} + \xi(t)\right) & \text{OFF state} \end{cases} \quad (1)$$

where $\xi(t)$ represents the random temperature fluctuations, which we generated as pseudorandom numbers having normal distribution. The equations for TCL dynamics have been broadly reviewed [20], [25]. When the air-conditioner is in the ON state, $\theta(t)$ decreases and TCL consumes constant power P=14 kW, as shown in Fig 1. When the air-conditioner is in OFF state, P=0.

The set point temperature is 20 $^0$C and the lower and upper dead band limits are 19.25 $^0$C and 20.75 $^0$C respectively, as shown in Fig 1. The ambient temperature is 32 $^0$C. R is the thermal resistance, and C is the thermal capacitance. Selecting R = 2 $^0$C/kW, C = 1.8 kWh/$^0$C, and P=14 kW, the cooling time is 20 minutes, and the heating time is 27 minutes. The total average TCL time period is 47 minutes. To model the load heterogeneity , the values of R and C for all the loads are calculated by adding pseudorandom values drawn from the standard uniform distribution on the open interval (0,1) to the aforementioned values [40].

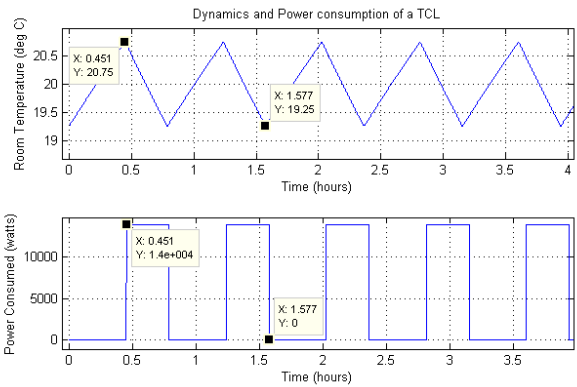

Fig 1. Dynamics of temperature (up) and power (down) consumption of a single TCL.

## III. UNSAFE SETPOINT CHANGE

At constant temperature, a heterogeneous population of TCLs achieves a steady state where instantaneous TCL states are uncorrelated and population consumes power with a constant rate up to small fluctuations due to the finite size of the TCL ensemble.

When used to make rapid changes in load, many straightforward methods for direct load control demonstrate large unwanted power oscillations due to correlations introduced by a synchronizing control signal [1], [41]. We will name control strategies leading to such oscillations as *unsafe*. One example of an unsafe method is control of a TCL ensemble by applying a sudden temperature set point change to all TCLs.

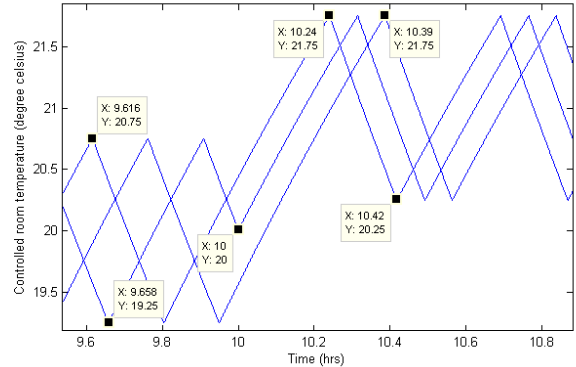

Fig. 2. Effect of the imposed upward shift in temperature set point on dynamics of temperatures of three TCLs.

Fig 2 shows the dynamics of temperatures for three TCLs to demonstrate the effect of a set point change. The initial set point temperature of 20 $^0$C is shifted to 21 $^0$C at time t = 10 hours. The dead band width is kept constant at 1.5 $^0$C, and the new dead band limits after the shift is applied are $[\theta_n, \theta_{n+}]$ = [20.25 $^0$C, 21.75 $^0$C]. After the control signal is applied (see Fig. 2), the two TCLs that were initially in the OFF state continue to remain in the OFF state until they reach the new upper dead band of 21.75 $^0$C after which their dynamics are constrained within new temperature dead band limits. The third TCL, which is initially ON, instantly turns OFF as its temperature at t =10 hours is less than 20.25 $^0$C.

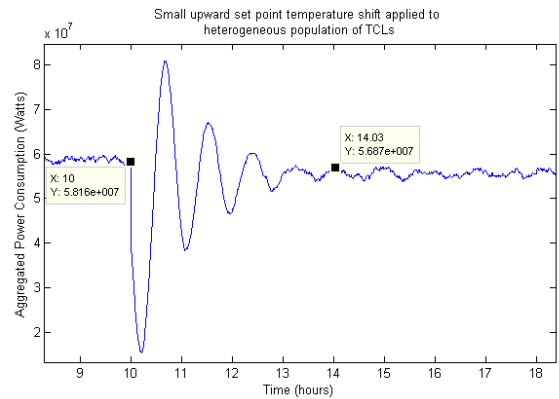

Fig 3. Aggregated response of power consumption of 10,000 TCLs to an upward shift by 1 $^0$C in temperature set point.

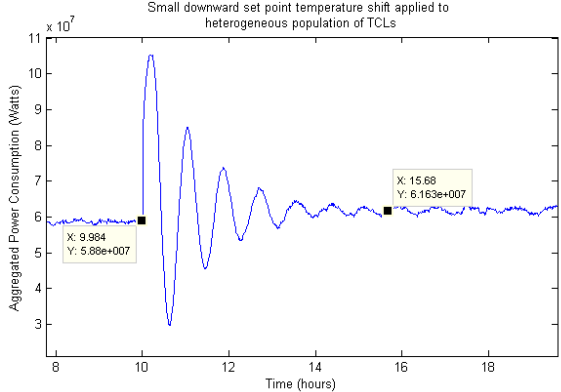

Fig 4. Aggregated response of power consumption of 10,000 TCLs to a downward shift by 1 °C in temperature set point.

Fig 3 and Fig 4 show the aggregated power consumed by a large ensemble of 10,000 TCLs, with parameters explained in Section II, when a set point shift of 1°C is suddenly applied in upward and downward directions, respectively. Large oscillations are observed for several hours due to synchronization of individual TCL states and the recovery of aggregated power to the steady state is sluggish.

## IV. GENERATION OF POWER PULSES USING SAFE PROTOCOL

We retain the simplicity of the open loop control methodology discussed above, but eliminate the problem of unwanted synchronization of TCLs by using timers installed in the temperature controllers. Using these timers, we are able to generate a power pulse and then slowly bring all TCLs back to their initial set points while avoiding unwanted oscillations. We refer to these methods as *timer-based safe protocols* where 'Safe' refers to the lack of oscillations. Here we will discuss two strategies, which we will name, respectively, SP-T1 and SP-T2, where SP stands for the "safe protocol" and T is for "timer".

### A. SP-T1 to delay the power consumption by TCLs

We start with SP-T1, a simpler strategy with a goal of delaying the power consumption of TCLs by a relatively large span of time (e.g. 1 hour). To achieve this, one does not have to switch TCLs to OFF for this duration of time. When SP-T1 is initiated by a single control signal all TCLs initially continue working as usual until they reach the point where they would normally switch from OFF to ON. At this point, the TCLs remain in the OFF state for an extra M minutes. Subsequently the TCLs should switch ON and return to their normal operation. Fig 5 illustrates the behavior with 3 TCLs for M = 10 minutes and the SP-T1 is initiated at t = 10 hours.

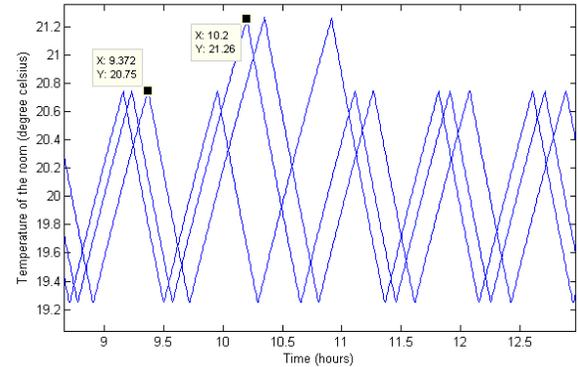

Fig 5. SP-T1 applied to 3 TCLs where TCLs are commanded to remain OFF for an extra M=10minutes

Conversely, if the goal is initially to absorb extra energy from the grid and then release it after a delay of ~30 minutes, TCLs should work as usual until they reach the point at which they would normally switch from ON to OFF. Instead, they continue in the ON state for M minutes before returning to the initially set parameters. Fig 6 illustrates the behavior with 3 TCLs for M = 10 minutes and the SP-T1 is initiated at t = 10 hours.

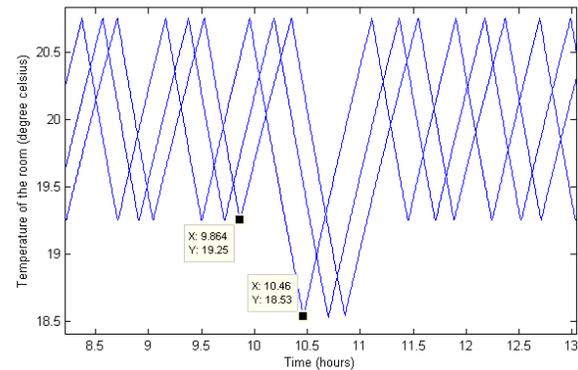

Fig 6. SP-T1 applied to 3 TCLs where TCLs are commanded to remain ON for an extra M=10 minutes

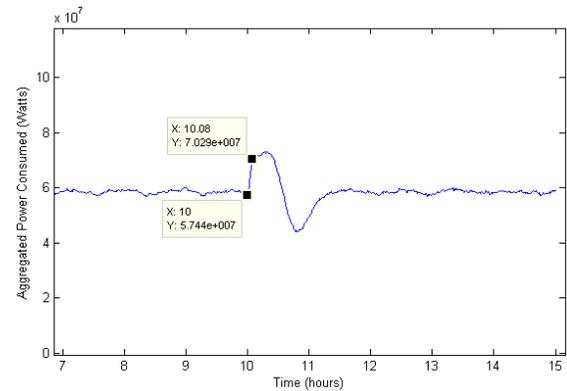

(a)

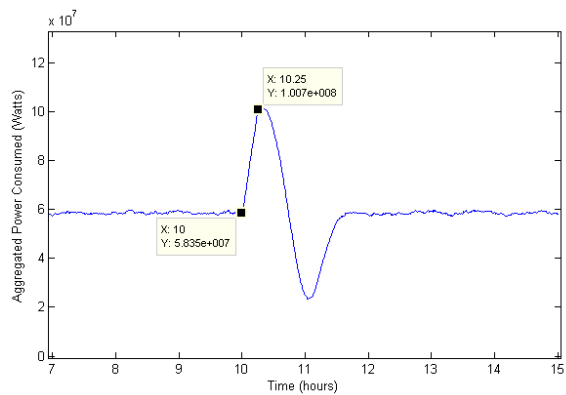

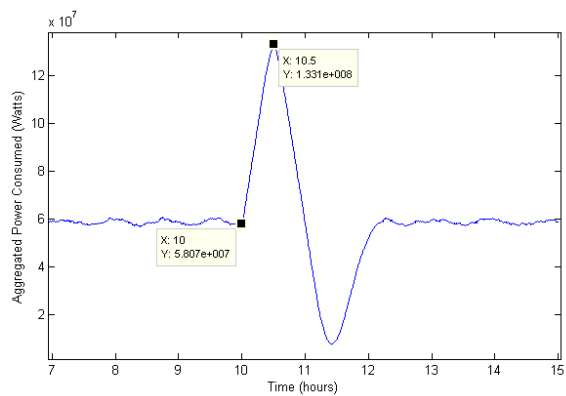

Fig 7. Aggregated power response of 10,000 TCLs. (a) M=5 minutes. (b) M=15 minutes. (c) M=30 minutes.

Fig 7a, Fig 7b and Fig 7c show the aggregated power for a heterogeneous population of 10,000 loads when SP-T1 is applied for M=5, 15 and 30 minutes, respectively. The power consumed increases linearly during the first M minutes as the TCLs remain in the ON state. If M>30min, all the loads finally appear in the ON state and the aggregated power consumed is close to the peak value of 140 MW (see Fig 7c). After achieving this maximum, the TCLs power consumption decreases almost linearly due to continuous return of TCLs to their original operation state. Eventually, the distribution of TCLs returns to the uncorrelated initial distribution, without producing any substantial power oscillations.

Here we note that the ensemble of TCLs behaves much more like batteries than standard power generators. Like batteries, for the TCLs to return to their original state, the energy that was initially absorbed (generated) by TCLs must be returned to (absorbed from) the grid. Inspection of Fig 5a, Fig 5b and fig 5c show that the shape and duration of these pulses are determined by the natural evolution of the TCL ensemble.

SP-T1 has advantages relative to the less sophisticated methods of TCL control. Comparing to Fig 3 where TCL set points is simply shifted, the additional timing input in SP-T1 contain the response to a well-defined time frame and eliminates the extended oscillations and the risk associated with these oscillations. The reason for the absence of remaining power oscillations in SP-T1 protocol is the same as for the safe protocol SP-2 described in [1].

The power pulses provided by SP-T1 are potentially useful for peak shaving or spinning reserve applications. However, the timescale of the response may better fit following power supply fluctuations, e.g. when power is provided by intermittent renewable sources.

### B. SP-T2 Protocol to generate a sharp power pulse

SP-T2 can be used to generate an abrupt power spike of precise duration and either sign. SP-T2 consists of the following steps:

Step 1: The SP-T2 control signal instructs the TCLs in the ON state to turn OFF for a prescribed time Δt (e.g. 2 to 3 minutes in our simulations). The TCLs also store the upper and lower dead band limits. There is a sudden drop in the aggregated power and all TCLs continue to remain OFF for time Δt.

Step 2: After time Δt, the TCLs that were turned OFF in Step 1 are instructed to turn ON and operate according to the original dead band limits. There will be a shift in the dead band position by the end of Step 2 for time Δt as shown in Fig 8. The size of this shift depends on the time for which the signal is applied and parameters of a TCL. The set point needs to be shifted back to the original value to ensure that the interference is not permanent. This is achieved in the next step.

Step 3: The TCLs that are commanded to turn ON in Step 2 are then made to operate according to the original dead band limits. However, the TCLs in OFF state should continue to remain OFF for time Δt after they reach the original upper dead band limit and only then they are instructed autonomously to operate according to the original dead band limits. Unlike Step 2, this step is slow ensuring that there is no subsequent sharp load spike and unwanted oscillations. After Step 3, all TCLs operate according to the original dead band limits

In order to better illustrate the dynamics of TCLs under this protocol, the time-dependence of the temperature of five TCLs is shown in Fig 8. At time t = 10 hours, signal is applied and all air-conditioners that are ON are commanded to switch OFF. Fig 8 shows that one TCL is ON when the signal is applied. This TCL remains OFF till time t = 10.04 hours (for around 2 minutes) and then turns ON and starts operating according to the original dead band limits.

Whereas, four TCLs that are OFF at time t = 10 hours continue to remain OFF for 2 minutes after they have reached the original upper dead band limit of 20.75 $^0$C, after which these TCLs turn ON and operate within original hysteresis band.

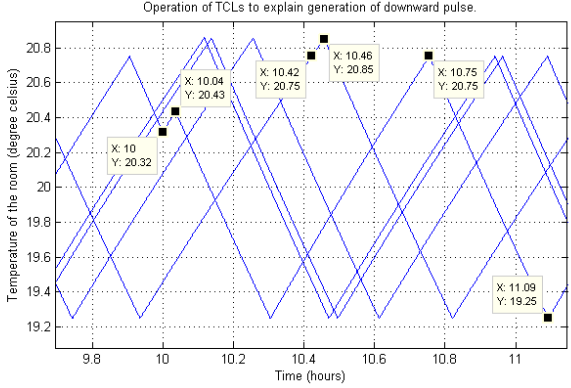

Fig 8. Temperature dynamics of five different TCLs under application of SP-T2 that generates a downward pulse.

It is observed that a set point shift of about 0.1 $^0$C occurs for the duration of the applied signal after which the set point is brought back to its original value. Thus, if the signal is applied for few minutes it will not cause the room temperature to change by a noticeable amount.

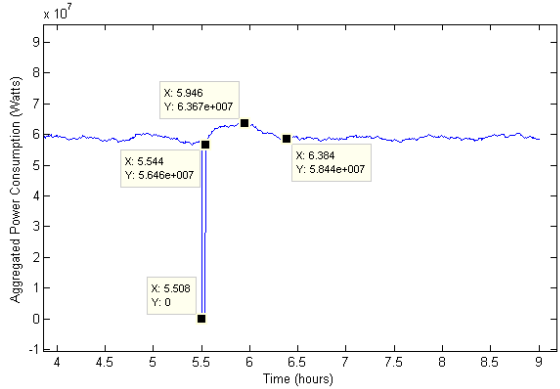

Fig 9. Downward pulse of 2 minutes duration generated by TCL population using SP-T2 safe protocol.

Fig 9 shows the simulation of a heterogeneous population of 10,000 TCLs. The power drops to zero when the signal is applied at time t = 5.508 hours. This stage continues for 2 minutes after which the aggregated power consumption increases to the initial value followed by slow return of TCLs to the steady state during time of a TCL cycle. Unwanted power oscillations are not observed and the peak amplitude, when the system tries to reach the equilibrium, is found to be a small fraction of the pulse generated. Similarly we can generate an upward pulse with duration of 2 minutes as shown in Fig 10 [1].

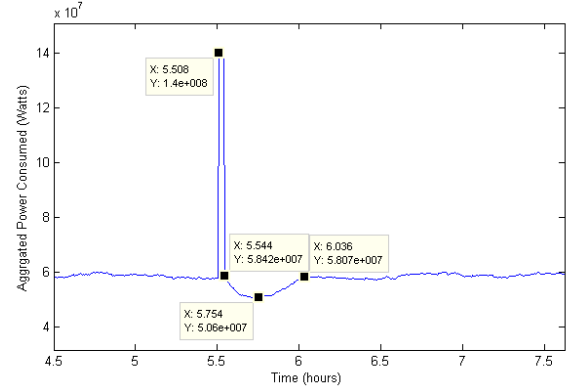

Fig 10. Upward pulse of 2 minutes duration generated by TCL population using SP-T2 safe protocol.

Here we note that, in reality, an instantaneous change of state of all TCLs cannot be achieved. While radio frequency signal can be transmitted almost instantly, TCLs may need time to safely change their mode. For example, for air-conditioners this time is typically about 20 seconds, which limits the duration and precision of SP-T2. Other electronic devices, without spinning mechanical components, such as water heaters may respond much faster. In what follows we assume such a TCL population having much smaller than 2 minutes time of response.

## V. OFFSETTING FAST FLUCTUATIONS

The safe protocol (SP-T2) from section IV has been briefly introduced in previous work [1], but its possible applications have not been discussed. In this section, we will look at how SP-T2 can help offset frequent fast fluctuations.

Consider a large heterogeneous population of TCLs. A group of TCLs are switched ON or OFF depending on the type (upward or downward) pulse to be generated. Having the knowledge of the net-demand profile from the utility, we apply control to different groups of TCLs thus avoiding frequent switching of the same set of TCLs. The magnitude of the pulse to be generated is directly proportional to the size of the ensemble. If $N$ is the number of TCLs in a group (we assume $N >> 1$). The power responses will depend linearly on $N$:

$$P_- = \alpha_- N, \qquad (2)$$

$$P_+ = \alpha_+ N, \qquad (3)$$

where $P_-$ and $P_+$ are magnitudes of respectively decrease and increase of power absorbed by TCLs when signal is sent to switch OFF and ON, respectively. Knowing N and the total power that ensemble consumes, constants $\alpha_-$ and $\alpha_+$ are proportionality coefficients, which can then be calculated.

We assume that the grid operator has access to such aggregate information to estimate $\alpha_+$ and $\alpha_-$.

We assume that the operator who sends control signals to TCLs has a short-term (e.g. 3 minutes) forecast for a size of a power fluctuation in the grid. Knowing the magnitude of the fluctuation that starts at time $t_1$, the number of TCLs (say $N_1$) to which the signal should be applied to command them to turn OFF can be calculated as

$$N_1 = \frac{P_{up}}{\alpha_-} \qquad (4)$$

where $P_{up}$ is the fluctuation in the upward direction. Now let at time $t_2$ the fluctuation occur in the downward direction. In order to offset it we will generate the pulse in the upward direction. Then, the number of TCLs (say $N_2$) that are commanded to remain ON for a short duration can be calculated as

$$N_2 = \frac{P_{dn}}{\alpha_+} \qquad (5)$$

where $P_{dn}$ is the fluctuation in the downward direction. In this way, control can be applied to different groups of TCLs until the sum of TCLs of all the groups equals the total number of TCLs in the ensemble. After this, the process can be restarted with TCLs that had time to return to the customer specified set points. Fig 11 and Fig 12 show the simulation results that demonstrate the effectiveness of SP-T2 in offsetting fast time scale fluctuations without leading to synchronization of individual states of TCLs.

*A) Heterogeneous Population of 15,000 TCLs*

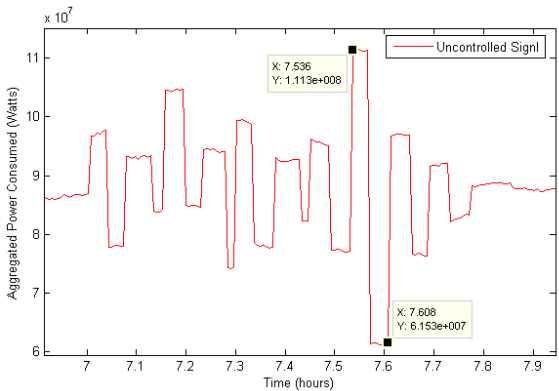

Fig 11. Aggregate power consumed by TCLs plus the step-like external fluctuations that operator desires to eliminate.

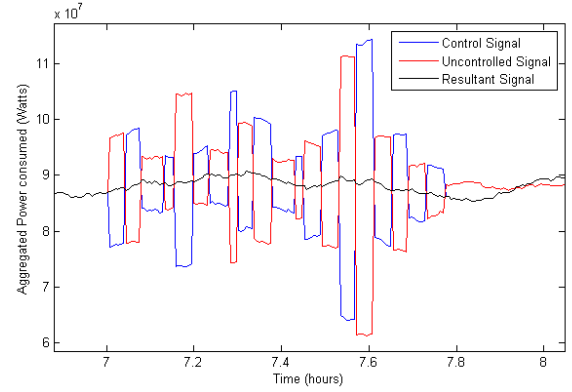

Fig 12. Offsetting fast time scale fluctuations. Blue curve is the power output of TCL ensemble that receives control signal to reduce fluctuations (Red). Black curve is the total power consumption, which illustrates the smoothing effect of control.

Fig 11 demonstrates results of our simulations of dynamics of a heterogeneous population of 15,000 TCLs. For t<7 hours, it shows the power demand of this TCL ensemble at the steady state (which would be approximately constant). For t>7 hours, a series of demand fluctuations occurs one after another. Then we assume that the control is applied from time t = 7 hours to time t = 7.8 hours. Fig 12 compares the power demand curve of Fig 11 (red) with the power consumed by a controlled population of TCLs (blue) and the total power demand generated (black). These results demonstrate that TCLs are able to respond quickly to offset fast fluctuations.

*B) Heterogeneous Population of 25,000 TCLs*

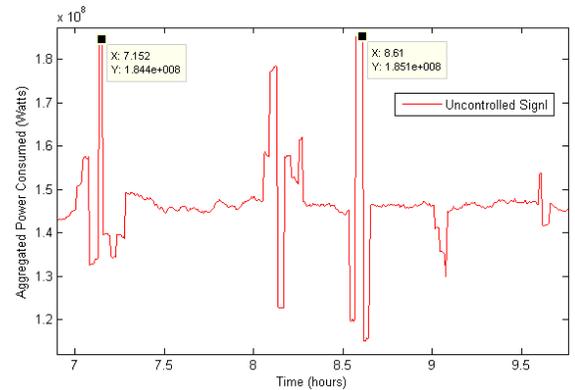

Fig 13. Aggregate power consumed by TCLs plus the arbitrary external fluctuations that operator desires to eliminate.

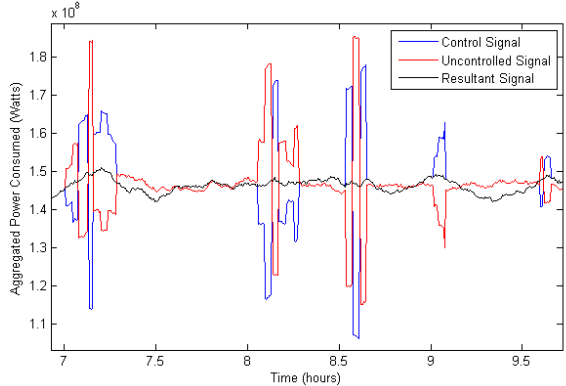

Fig 14. Offsetting fast time scale fluctuations. Blue curve is the power output of TCL ensemble that receives control signal to reduce fluctuations (Red). Black curve is the total power consumption, which illustrates the smoothing effect of control.

Figure 13 corresponds to a different case of intermittent fluctuations with random shapes. As size of the ensemble is bigger, the magnitude of a pulse that can be generated to offset the fluctuation is also large. In this example, pulses of magnitude 35 MW to 40 MW are generated. Figure 14 shows the versatility of SP-T2 lies in the fact that any arbitrary pulse can be generated for short duration without leading to any parasitic oscillations.

## VI. NON-EQUILIBRIUM CONDITIONS

In the previous sections, as well as in previous studies of safe protocols [1], [20], [25], the ambient temperature was considered constant and the response to the control signal was studied assuming that aggregated power consumption was initially fully equilibrated. In reality, ambient temperature continuously changes during the day. If the degree of randomness is not sufficiently strong, the assumption that the ensemble is near equilibrium may not apply. Here, we investigate the application of a safe protocol to generating power pulses while ambient outdoor temperature is changing with time. The TCL ensemble has strong heterogeneity, and the natural stochastic fluctuations, which normally are the main force for equilibrating the ensemble, were set to zero.

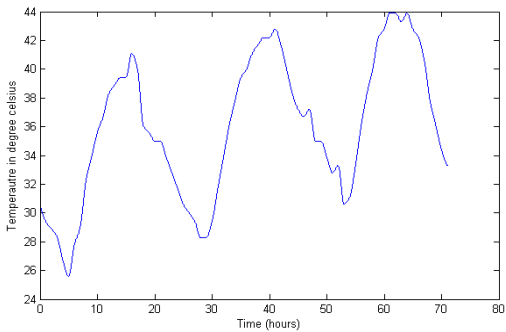

Fig 15. Ambient temperature variation in Arizona.

Fig 15 shows the ambient temperature variation in Arizona during summer 2012 for a three days period (72 hours) [42]. We simulated the operation of a heterogeneous population of 10,000 TCLs, assuming that they have been working continuously during the same period of time with fixed customer set points.

We again assume that the control authority has an estimate of the total number of working air-conditioners and the power they consume at every moment. At one moment in time (t = 50 hours), we assumed that the utility company needs to generate a short (2 minutes) but strong power pulse, (either with positive or negative sign). Our goal was to test proper functionality of safe protocol under such non-equilibrium conditions.

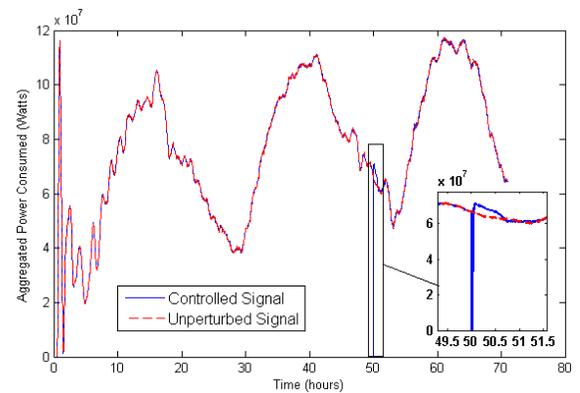

Fig 16. Aggregate power response of the controlled (blue) and the unperturbed (red) TCL ensemble when a downward pulse is generated.

Fig 16 shows our results for the aggregated power consumed by such an ensemble. The blue signal represents the controlled signal and the red signal is the signal at no control. In the blue trace, all TCLs are subjected to a 2 minute SP-T2 at time t = 50 hours to generate downward pulse and are then returned to the initial set point. It shows, in particular, that the broad heterogeneity of the ensemble, alone, is sufficient to make power consumption of TCLs relatively smooth while clearly following the trend of the outdoor temperature dynamics. SP-T2 works as desired, i.e. we did not observe any additional side effects on top of the natural power demand evolution and fluctuations of the ensemble. We observe that both signals trace the same path ensuring proper working of the safe protocol as described in section IV. The mean, standard deviation (SD) of distribution of different TCL parameters are tabulated in Table 1.

TABLE 1

Mean and standard deviation of different parameters when a downward pulse
is created.

| Upper dead band limit (UDL) $^0$C | | Lower dead band limit (LDL) $^0$C | | Set point (SP) $^0$C | | Dead band width (DW) $^0$C | |
|---|---|---|---|---|---|---|---|
| Mean | SD | Mean | SD | Mean | SD | Mean | SD |
| 21.2463 | 0.2895 | 19.2435 | 0.4076 | 20.2449 | 0.3230 | 2.0027 | 0.2873 |

Figure 17 shows the similar response to the upward pulse with similar conclusions. The mean, standard deviation (SD) of different parameters are tabulated in Table 2.

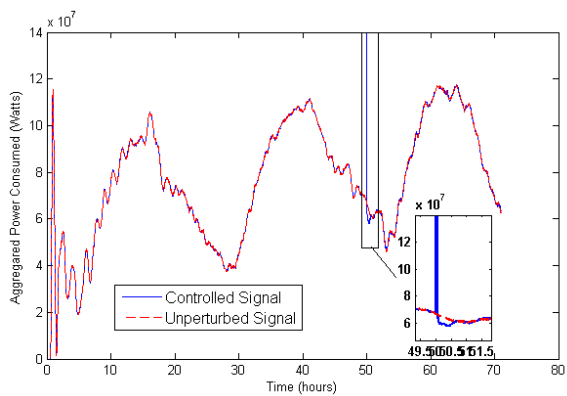

Fig 17. Aggregate power response of the controlled (blue) and the unperturbed (red) TCL ensemble when a downward pulse is generated.

TABLE 2

Mean and standard deviation of different parameters when a upward pulse is
created.

| Upper dead band limit (UDL) $^0$C | | Lower dead band limit (LDL) $^0$C | | Set point (SP) $^0$C | | Dead band width (DW) $^0$C | |
|---|---|---|---|---|---|---|---|
| Mean | SD | Mean | SD | Mean | SD | Mean | SD |
| 21.2521 | 0.2891 | 19.2552 | 0.4052 | 20.2537 | 0.3207 | 1.9969 | 0.2899 |

## VII. CONCLUSION

In this work, we showed that adding timers to the endpoint temperature controller can be used to generate upward or downward pulses of large magnitude and small duration making it suitable to offset fast time scale fluctuations with knowledge of only several basic aggregate parameters that describe the state of the ensemble. Pulses of varying shapes, magnitudes and duration can be generated with no synchronization of individual states of TCLs and associated power oscillations at the end of the protocol. We also demonstrated that even during considerable variation of the outdoor temperature, safe protocols produce the desired response without any visible unwanted effects.

The TCL control protocols described in this manuscript are minimally invasive and less expensive to implement. For example, they only require one-way communication and do not require the collection of data for each individual TCL. Also, the set point can be guaranteed to change within a comfortable level when the control signal is applied with the set points returning to their original value after the needed power pulse is generated. Such a minimally invasive strategy will avoid customer dissatisfaction. We expect that strategies based on the safe protocols should become commercially interesting when the number of TCLs in a community under control exceeds several thousands.

Our results suggest that it can be economically attractive for utility companies to develop demand side management programs in which few basic upgrades, such as timers and a small amount of memory, are imbedded in TCLs to instruct them to turn ON or OFF according to safe protocol strategies. Such demand side control strategies will help to reduce cost and maintain grid reliability.

## REFERENCES

[1] N. Sinitsyn, S. Kundu, and S. Backhaus, "Safe protocols for generating power pulses with heterogeneous populations of thermostatically controlled loads," *Energy Conversion and Management*, 2013.

[2] M. Klobasa, "Analysis of demand response and wind integration in Germany's electricity market," *IET Renewable Power Generation*, vol. 4, no. 1, p. 55, 2010.

[3] R. M. Delgado, "Demand-side management alternatives," *Proceedings of the IEEE*, vol. 73, no. 10, pp. 1471–1488, 1985.

[4] Z. Xu, J. Ostergaard, and M. Togeby, "Demand as frequency controlled reserve," *Power Systems, IEEE Transactions on*, vol. 26, no. 3, pp. 1062–1071, 2011.

[5] W. Burke and D. Auslander, "Robust control of residential demand response network with low bandwidth input," *ASME Dynamic Systems and Control Conference*, pp. 3–5, 2008.

[6] N. Ruiz, I. Cobelo, and J. Oyarzabal, "A Direct Load Control Model for Virtual Power Plant Management," *IEEE Transactions on Power Systems*, vol. 24, no. 2, pp. 959–966, May 2009.

[7] A. Molina-García and M. Kessler, "Probabilistic characterization of thermostatically controlled loads to model the impact of demand response programs," *Power Systems, IEEE Transactions on*, vol. 26, no. 1, pp. 241–251, 2011.

[8] S. El-Ferik, "Identification of physically based models of residential air-conditioners for direct load control management," *Control Conference, 2004. 5$^{th}$ Asian*, 2004.

[9] G. Heffner, C. Goldman, B. Kirby, and M. Kintner-Meyer, "Loads providing ancillary services: Review of international experience," vol. 11231, no. May, 2008.

[10] G. Strbac, "Demand side management: Benefits and challenges," *Energy Policy*, vol. 36, no. 12, pp. 4419–4426, Dec. 2008.

[11] D. Kirschen and G. Strbac, *Fundamentals of Power System Economics*. Chichester, UK: John Wiley & Sons, Ltd, 2004.

[12] K. Schisler, T. Sick, and K. Brief, "The role of demand response in ancillary services markets," *Transmission and Distribution Conference and Exposition*, pp. 1–3, 2008.

[13] C. Sastry, R. Pratt, V. Srivastava, and S. Li, *Use of Residential Smart Appliances for Peak-Load Shifting and Spinning Reserves: Cost/Benefit Analysis*, no. December. 2010.

[14] J. Laurent, G. Desaulniers, R. Malhame, and F. Soumis, "A column generation method for optimal load management via control of electric water heaters," *Power Systems, IEEE Transactions on*, 1995.

[15] N. Lu and S. Katipamula, "Control Strategies of Thermostatically Controlled Appliances in a Competitive Electricity Market," *IEEE Power Engineering Society General Meeting, 2005*, pp. 164–169.

[16] D. S. Callaway and I. A. Hiskens, "Achieving Controllability of Electric Loads," *Proceedings of the IEEE*, vol. 99, no. 1, pp. 184–199, Jan. 2011.

[17] D. S. Callaway, "Tapping the energy storage potential in electric loads to deliver load following and regulation, with application to wind energy," *Energy Conversion and Management*, vol. 50, no. 5, pp. 1389–1400, May 2009.

[18] K. Kalsi, M. Elizondo, J. Fuller, S. Lu, and D. Chassin, "Development and Validation of Aggregated Models for Thermostatic Controlled Loads with Demand Response," *2012 45th Hawaii International Conference on System Sciences*, pp. 1959–1966, Jan. 2012.

[19] R. Malhame, "Electric load model synthesis by diffusion approximation of a high-order hybrid-state stochastic system," *IEEE Transactions on Automatic Control*, vol. 30, no. 9, pp. 854–860, Sep. 1985.

[20] S. Kundu and N. Sinitsyn, "Safe protocol for controlling power consumption by a heterogeneous population of loads," *American Control Conference (ACC), Montreal, 2012*, pp. 2947–2952, 2012.

[21] J. L. Mathieu and D. S. Callaway, "State Estimation and Control of Heterogeneous Thermostatically Controlled Loads for Load Following," *2012 45th Hawaii International Conference on System Sciences*, pp. 2002–2011, Jan. 2012.

[22] C. Perfumo, E. Kofman, J. H. Braslavsky, and J. K. Ward, "Load management: Model-based control of aggregate power for populations of thermostatically controlled Loads," *Energy Conversion and Management*, vol. 55, pp. 36–48, Mar. 2012.

[23] C. Ucak and R. Caglar, "The Effects of Load Parameter Dispersion and Direct Load Control Actions on Aggregated Load," *Power System Technology, 1998. Proceedings*, pp. 280–284, 1998.

[24] S. Bashash and H. Fathy, "Modeling and control insights into demand-side energy management through setpoint control of thermostatic loads," *American Control Conference (ACC), 2011*, pp. 4546–4553, 2011.

[25] S. Kundu, N. Sinitsyn, S. Backhaus, and I. Hiskens, "Modeling and control of thermostatically controlled loads," *arXiv preprint arXiv: 1101.2157*, 2011.

[26] E. Bompard, E. Carpaneto, G. Chicco, and R. Napoli, "Analysis and Modelling of Thermostatically-Controlled Loads," *Electrotechnical Conference , 1996. MELECON'96., 8$^{th}$ Mediterranean*, no. 39, 1996.

[27] D. Chassin and J. Malard, "The Equilibrium Dynamics of Thermostatic End-use Load Diversity as a Function of Demand," *arXiv preprint nlin/0409037*, pp. 1–9, 2004.

[28] W. Zhang, K. Kalsi, J. Fuller, M. Elizondo, and D. Chassin, "Aggregate Model for Heterogeneous Thermostatically Controlled Loads with Demand Response," *Power and Energy Society General Meeting*, pp. 1–8, 2012.

[29] C. Chong and A. Debs, "Statistical synthesis of power system functional load models," *Decision and Control including the Symposium on Adaptive Processes*, 1979.

[30] S. Koch, J. Mathieu, and D. Callaway, "Modeling and control of aggregated heterogeneous thermostatically controlled loads for ancillary services," *Proc. PSCC*, 2011.

[31] N. Lu, D. Chassin, and S. Widergren, "Modeling Uncertainties in Aggregated Thermostatically Controlled Loads Using a State Queueing Model," *Power Systems, IEEE Transactions on*, 2005.

[32] R. Mortensen and K. Haggerty, "A Stochastic Computer Model for Heating and Cooling Loads," *Power Systems, IEEE Transaction on*, vol. 3, no. 3, 1988.

[33] K. Kalsi, F. Chassin, and D. Chassin, "Aggregated Modeling of Thermostatic Loads in Demand Response: A Systems and Control Perspective," *Decision and Control and European Control Conference (CDC-ECC)*, 2011.

[34] J. Bendtsen and S. Sridharan, "Efficient Desynchronization of Thermostatically Controlled Loads," *arXiv preprint arXiv:1302.2384*, 2013.

[35]     "U.S. Energy Information Administration." [Online].
         Available: http://www.eia.gov/.

[36]     "Cooling a Warming Planet: A Global Air
         Conditioning Surge," 2012. [Online]. Available:
         http://e360.yale.edu/feature/cooling_a_warming_plan
         et_a_global_air_conditioning_surge/2550/.

[37]     B. Kirby, *Spinning reserve from responsive loads*, no.
         March. 2003.

[38]     B. Kirby, J. Kueck, T. Laughner, and K. Morris,
         "Spinning Reserve from Hotel Load Response," *The
         Electricity Journal*, vol. 21, no. 10, pp. 59–66, Dec.
         2008.

[39]     B. Kirby, "Load response fundamentally matches
         power system reliability requirements," *Power
         Engineering Society General Meeting, 2007. IEEE*,
         pp. 1–6, 2007.

[40]     "MATLAB." The MathWorks Inc., Natick,
         Massachusetts, 2011.

[41]     S. Ihara and F. Schweppe, "Physically based
         modeling of cold load pickup," *Power Apparatus and
         Systems, IEEE Transactions on*, no. 9, pp. 4142–
         4150, 1981.

[42]     "Weather Underground." [Online]. Available:
         http://www.wunderground.com/.